
\input amstex
\documentstyle{amsppt}
\NoBlackBoxes
\NoRunningHeads


\define\ggm{{G}/\Gamma}

\define\vre{\varepsilon}
\define\hd{Hausdorff dimension}
\define\hs{homogeneous space}
\define\df{\overset\text{def}\to=}
\define\un#1#2{\underset\text{#1}\to#2}
\define\br{\Bbb R}
\define\bn{\Bbb N}
\define\bz{\Bbb Z}

\define\di{Diophantine}
\define\da{Diophantine approximation}

\define\vv{\bold v}

\define\vw{\bold w}

\define\vt{\bold t}
\define\vu{\bold u}
\define\g{\frak g}
\define\h{\frak h}
\define\fa{\frak a}

\define\Lie{\operatorname{Lie}}
\define\Lip{\operatorname{Lip}}
\define\cic{C^\infty_{comp}}


\define\mn{{m+n}}
\define\nz{\smallsetminus \{0\}}

\define\cag{$(C,\alpha)$-good}

\define\ssm{\smallsetminus}

\define\supp{\operatorname{supp}}
\define\const{\operatorname{const}}
\define\dist{\operatorname{dist}}
\define\SL{\operatorname{SL}}
\define\GL{\operatorname{GL}}
\define\Ad{\operatorname{Ad}}
\define\mr{M_{m,n}}
\define\amr{$A\in M_{m,n}(\br)$}

\define\ehs{expanding horospherical subgroup}

\define\diag{{\operatorname{diag}}}
\define\vo{{\operatorname{vol}}}
\newif\ifdraft\drafttrue


\font\sb = cmbx8 scaled \magstep0

\long\def\comdima#1{\ifdraft{\sb #1 }\else\ignorespaces\fi}

\topmatter
\title On effective equidistribution
of expanding translates of certain orbits in the space of lattices
 \endtitle 

\author { D.\,Y. Kleinbock and G.\,A. Margulis} \endauthor



    \address{ Department of
Mathematics, Brandeis University, Waltham, MA 02454}
  \endaddress

\email kleinboc\@brandeis.edu \endemail

    \address{ Department of Mathematics, Yale University, 
   New Haven, CT 06520}
  \endaddress

\email margulis\@math.yale.edu \endemail
           
\abstract 
We prove an effective version of a result obtained earlier by Kleinbock and Weiss \cite{KW} on equidistribution of expanding translates of orbits of horospherical subgroups in the space of lattices.
\endabstract


\endsubjclass

\endtopmatter

\document

\heading{1.  Introduction}\endheading

The motivation for this work is a result obtained recently in \cite{KW}.
Fix $m,n\in\bn$, set $k = \mn$ and let
$$
G = \SL_{k}(\br),\  \Gamma = \SL_{k}(\bz),\ u_{ \sssize Y} = \pmatrix
I_{m} & Y
\\ 0 & I_{n}\endpmatrix,\ 
H = \left\{\left.u_{ \sssize Y}\right|Y\in\mr\right\},\tag 1.1
$$ where $\mr$ stands for 
the space of $m\times n$ matrices with real entries. 
Then $H$ is a unipotent abelian subgroup of $G$ which is
{\sl expanding horospherical\/} with respect to 
$$g_t =  \diag(e^{t/m}, \ldots,
e^{t/m}, e^{-t/n}, \ldots,
e^{t/n})\,,\quad t > 0\,.\tag 1.2
$$
The latter, by definition, means that the Lie algebra 
of $H$ is 
the span of eigenspaces of $\Ad(g_t)$, $t > 0$,
with eigenvalues bigger than $1$ in absolute value.

The space $X\df\ggm$ can be identified with the space of unimodular
lattices in $\br^k$, on which $G$ acts by left translations.
Denote by $\pi$ the natural projection $G\to X,\ g \mapsto g\Gamma$, and
for any $z\in X$ let $\pi_z:G\to X$ be defined by $\pi_z(g) = gz$.
Also denote by $\bar \mu$ the $G$-invariant probability measure on $X$
and by $\mu$ the Haar measure on $G$ such that $\pi_*\mu = \bar \mu$.
Fix a Haar measure $\nu$ on $H$. Note that the $H$-orbit foliation
is unstable with respect to the action of $g_t$, $t > 0$.
It is well known 
that  for any Borel
probability  measure $\nu'$ on $H$ absolutely
continuous with respect 
to $\nu$ and for any $z\in X$, $g_t$-translates of  $(\pi_z)_*\nu'$ 
become equidistributed, that is, weak-$*$ converge to 
$\bar\mu$ as $t
\to\infty$. An effective version of this 
statement was obtained in \cite{KM1, Proposition 2.4.8}. In order to state
that result, it will be convenient to introduce the
following notation: for  $f\in L^1(H,\nu)$, a bounded continuous 
function $\psi$
on
$X$, 
$z\in X$ and $g\in G$ define
$$I_{f,\psi}(g,z)\df \int_H
f(h)\psi(g_thz)\,d\nu(h)\,.$$
In other words, $I_{f,\psi}(g,z)$ is the result of evaluation of
the $g$-translate of $(\pi_z)_*\nu'$ at  $\psi$, where
$d\nu' = f\, d\nu$. Then equidistribution of $g_t$-translates of
$(\pi_z)_*\nu'$ amounts  to the convergence of $I_{f,\psi}(g_t,z)$
to $\int_H
f\cdot\int_X\psi
$ as $t\to\infty$ (unless it causes confusion, we will  omit
measures in the integration
 notation for the sake of brevity).

\medskip
The following is 
a slightly
simplified
form of 
\cite{KM1, Proposition 2.4.8}:

\proclaim{Theorem 1.1} There exists $\,\gamma > 0$ such that
for any $f\in \cic(H)$, 
$\psi\in
\cic(X)$
 and for any compact subset $L$ of $X$ there
exists  a constant $C=C(f,\psi,L)$ such that  for
all $z\in L$ and  any $\,t\ge 0$
$$
\left|I_{f,\psi}(g_t,z) - \int_H
f\,
\int_X\psi\,
\right|\le C e^{-\gamma
t}\,.\tag 1.3  
$$
\endproclaim

The proof used the exponential decay of correlations of the $G$-action
on $X$ (called `condition (EM)' in \cite{KM1}). See \S 2 for more detail.

\medskip

Motivated by some questions in simultaneous \da, the first named
author and Barak Weiss considered translates of $H$-orbits on $X$
by diagonal elements of $G$ other than $g_t$. 
Specifically, following \cite{KW}, let us denote by $\fa^+$
the set of $k$-tuples $\vt = (t_1,\dots,t_{k})\in \br^{k}$
such that
$$
t_1,\dots,t_{k} > 0
\quad \text{and}\quad 
\sum_{i = 1}^m t_i =\sum_{j = 1}^{n} t_{m+j} \,,
$$ 
and for $\vt\in\fa_+$ define
$$g_\vt \df \diag(e^{t_1}, \ldots,
e^{ t_m}, e^{-t_{m+1}}, \ldots,
e^{-t_{k}})\in G\tag 1.4$$
and 
$$
\lfloor \vt\rfloor \df \min_{i = 1,\dots,k} t_i 
$$ 
(the latter, roughly speaking, measures the distance between $\vt$
and the walls of the cone  $\fa^+\subset \br^{k}$).

The theorem below is a reformulation of \cite{KW, Theorem 2.2}:

\proclaim{Theorem 1.2} For any $f\in L^1(H,\nu)$, any
continuous compactly supported $\psi:X\to\br$, any compact subset $L$ of
$X$
 and any $\,\vre>0$ there exists $T > 0$ such that 
$$
\left|
I_{f,\psi}(g_\vt,z)
- \int_H
f\,
\int_X\psi\,
\right| < \vre\,
$$
for all ${z}\in L$ and $\vt\in\fa^+, \, \lfloor \vt\rfloor\ge T$.
\endproclaim

That is, $g_{\vt}$-translates of $H$-orbits become
equidistributed as
$\lfloor
 \vt\rfloor
\to\infty$ uniformly in ${z}$ when the latter
is restricted to compact subsets of $X$. The proof
relies on S.\,G.~Dani's
classification of measures invariant under horospherical subgroups and
the so-called `linearization method'.
The purpose of the present paper is to prove an effective
version of the above theorem:

\proclaim{Theorem 1.3} There exists $\,\tilde \gamma > 0$ such that
for any $f\in \cic(H)$, 
$\psi\in
\cic(X)$
 and for any compact  $L\subset X$ there
exists   $\tilde C=\tilde C(f,\psi,L)$ such that  for
all $z\in L$ and all $\vt\in\fa^+$
$$
\left|
I_{f,\psi}(g_\vt,z) - \int_H
f\,
 \int_X\psi\,
\right|\le \tilde C e^{-\tilde\gamma
\lfloor
 \vt\rfloor}\,.
$$
\endproclaim

Note that the above statement follows from Theorem 1.1 when $k = 2$, that is, $G = \SL_2(\br)$,
but is new for $k > 2$.
The proof uses the `exponential mixing' approach of \cite{KM1, KM3} together
with effective nondivergence estimates obtained in \cite{KM2}.
We will describe these two parts in \S\S\ 2 and 3 respectively, and
then 
proceed with  the proof of Theorem 1.3 in
\S 4. We remark that the method of proof readily extends 
to the set-up 
more general than (1.1). Note also that, as observed by N.\ Shah in \cite{S, Remark 1.0.2},  Theorem 1.3 can be used to strengthen one of the main results of \cite{KW}, that is,
 \cite{KW, Theorem 1.4}, which constitutes a \di\ application of Theorem 1.2.

\medskip

 {\bf Acknowledgements.} The authors are grateful to the Fields Institute  for Research in Mathematical Sciences (Toronto, Canada),
 where this project has commenced, and to the referee for useful remarks.
 The work of the first named author was 
supported in part by NSF
Grants DMS-0239463 and DMS-0801064,
and that of the second author by NSF Grants DMS-0244406 and DMS-0801195.

\heading{2.  Exponential mixing and $g_t$-translates}\endheading

{\bf Notation:} We will fix a
right-invariant metric `dist' on $G$, giving rise to the corresponding
metric on $X$.
$B(x,r)$ will stand for an open ball of radius $r$ centered at $x$.
If a metric space is $G$ or its subgroups, we will abbreviate
$B(e,r)$ to $B(r)$. When
necessary, we will use subscripts denoting the ambient metric spaces.
$\|\cdot\|_\ell$ and $\|\cdot\|_{C^\ell}$ will stand for the $(2,\ell)$-Sobolev and $C^l$ norms
respectively. 
We define 
$$W^{2,\infty}(X) = \{\psi \in C^\infty(X): \|\psi\|_\ell < \infty\ \forall\,\ell\in\bn\}\,;
$$
clearly  $
\cic(X)\subset W^{2,\infty}(X)$.
In fact, $
W^{2,\infty}(X)$ coincides with the set of functions $\psi\in C^\infty(X)$
such that $D\psi\in L^2(X)$ for any
$D$ from the universal  enveloping  algebra of  $\Lie(G)$. We  let
$\langle\cdot,\cdot\rangle$ stand for the inner product in
$L^2(X)$.  We also denote by $\|\psi\|_{\Lip}$ the Lipschitz constant of 
a function $\psi$ on $X$, 
$$\|\psi\|_{\Lip}\df \sup_{x , y\in X,\ x\ne y}\frac{|\psi(x) -
\psi(y)|}{\dist(x,y)}\,,$$ and let 
$\Lip(X)\df \{\psi: \|\psi\|_{\Lip} < \infty\}$.

The following property of the $G$-action on $X$ is deduced in \cite{KM3}
from the 
spectral gap 
on $L^2(X)$: 

\proclaim{Theorem 2.1 { \rm[KM3, Corollary 3.5]}} There
exist  $\gamma>0$ and $\ell\in\Bbb N$ such that for
any two functions $\varphi,\psi\in 
W^{2,\infty}(X)$ 
and for any $\,t\ge 0$ one has 
$$
\left|\langle g \varphi,\psi\rangle -  \int_X \varphi\, 
\int_X\psi\, 
\right| 
\ll \|\varphi\|_{\ell}\|\psi\|_{\ell}\cdot e^{-\gamma \dist(g,e)}\,. 
$$
\endproclaim

Here and hereafter
the 
implicit constants in
$\ll$  depend only on the dimensions of the corresponding spaces
and the choices of the metric. Taking $g = g_t$ as in (1.2), it follows that
$$
\left|\langle g_t\varphi,\psi\rangle -  \int_X \varphi\, 
\int_X\psi\, 
\right| 
\ll \|\varphi\|_{\ell}\|\psi\|_{\ell}\cdot e^{-\gamma t}\,. \tag
2.1  
$$

\medskip 
An estimate analogous to (2.1) was used in  \cite{KM1} to derive Theorem 1.1.
In this section we apply Theorem 2.1 to prove a statement similar to
 Theorem 1.1,  providing some information as to how
$C$ in (1.3) 
depends on $f$ and $L$. 
The argument follows the lines of the proof in \cite{KM1}; in fact, the
statement below is basically an intermediate step in the proof of
\cite{KM1, Proposition 2.4.8}. However we have decided to include details
for the sake of making this paper self-contained.

\medskip

To pass from $\langle g_t\varphi,\psi\rangle$ to 
$I_{f,\psi}(g_t,z)$, we need to thicken $f$ into a 
suitable function $\varphi$ on $X$. To explain this process, we need to
introduce some more notation. Let
$$
H^- = \left\{\left. \pmatrix
I_{m} & 0
\\ Y & I_{n}\endpmatrix\right|Y\in M_{n,m}\right\}
$$
and
$$
H^0 = \left\{\left. \pmatrix
A & 0
\\ 0 & B\endpmatrix\right|A\in \GL_m(\br), \, B\in \GL_n(\br),
\,\det(A)\det(B) = 1\right\}\,.
$$
The product map $H^- \times H^0 \times H\to G$ is a local diffeomorphism;
we will choose $r_0$ so that the inverse of this map is 
well defined on $B_G(r_0)$.
Note that $H^-$ is expanding horospherical with respect to $g_{-t}$, $t
> 0$, while  $H^0$ is centralized by $\{g_t\}$. Thus, the inner
automorphism
$\Phi_t$ of $G$ given by  $\Phi_t(g)\df g_t h (g_t)^{-1}$ is 
non-expanding
on  the group 
$$
\tilde H \df H^- 
H^0 = \left\{\left. \pmatrix
A & 0
\\ Y & B\endpmatrix\right.
\right\}\,;
$$
in fact, one has 
$$ 
\forall\,r > 0\ \forall\,t>0\quad
\Phi_t\big(B_{\tilde H}({r})\big)\subset
B_{\tilde H}(r)\,.\tag 2.2
$$

Let us choose Haar measures $\nu^-$, $\nu^0$ on  $H^-$, $H^0$ 
respectively,  normalized  so that 
$\mu$ is 
locally almost the product of $\nu^-$, $\nu^0$ and $\nu$.
 By the
latter, in view of \cite{B,
Ch.~VII, \S 9, Proposition 13},
we mean that $\mu$ can be expressed via $\nu^-$, $\nu^0$ and $\nu$ in the
following way: for any $\varphi\in L^1(G)$
$$
{\int_{H^-  H^0  H} \varphi(g)\, d\mu(g)} = 
{\int_{H^-\times H^0\times H}\varphi(h^-h^0h)\Delta(h^0)\,d\nu^-(h^-)\,
d\nu^0(h^0) \, d\nu(h)}\,,\tag 2.3  
$$
where $\Delta$ is the modular function of (the non-unimodular group)
$\tilde H$. 

\medskip

The `thickening' will be based on the following properties of the Sobolev
norm, cf.\
\cite{KM1, Lemma 2.4.7}:

\proclaim{Lemma 2.2} {\rm (a)} 
For any
$\,r > 0$, 
there exists a nonnegative function $\theta\in
\cic(\br^N)$ such that 
{\rm supp}$(\theta)$ is inside $B(r)$, $\int_{\br^N} \theta
= 1$, and
$\|\theta\|_{\ell}\ll   r^{-(\ell+N/2)}$. 

{\rm (b)} Given 
$\theta_1\in
\cic(\br^{N})$, $\theta_2\in\cic(\br^{N})$, define
$\theta\in
\cic(\br^{N})$ by \newline $\theta(x) =
\theta_1(x)\theta_2(x)$. Then
$\|\theta\|_{\ell}\ll
\|\theta_1\|_{\ell}
\|\theta_2\|_{C^\ell}$.  

{\rm (c)} Given 
$\theta_1\in
\cic(\br^{N_1})$, $\theta_2\in\cic(\br^{N_2})$, define
$\theta\in
\cic(\br^{N_1+N_2})$ by $\theta(x_1,x_2) =
\theta_1(x_1)\theta_2(x_2)$. Then
$\|\theta\|_{\ell}\ll
\|\theta_1\|_{\ell}
\|\theta_2\|_{\ell}$.  
\endproclaim

We will
apply
the above lemma to functions supported on small enough balls centered at
identity elements in $G$, $H$, $H^0$,  $H^-$.

\proclaim{Theorem 2.3} 
Let $f\in \cic(H)$,
$0 < r < r_0/2$ and $z\in X$ be such that

{\rm (i)} $\supp f \subset B_H(r)$, and

{\rm (ii)} $\pi_z$ is injective on $B_G(2r)$.

\noindent Then for any 
$\psi\in
W^{2,\infty}(X)\cap \Lip(X)$ with $\int_X\psi = 0$ there exists $E = E(\psi)$ such that for 
 any $\,t\ge 0$ one has
$$
\left| I_{f,\psi}(g_t,z)
\right|\le\ E \left(r 
\int_H |f|\,
 +
r^{-(2\ell + N/2)} \|f\|_{\ell}
e^{-\gamma
t}\right)\,,
\tag 2.4  
$$
where $\gamma$ and $\ell$ are as in Theorem 2.1 and $N = m^2 + mn + n^2
-1 = \dim \tilde H$. 
\endproclaim


\demo{Proof} 
Using Lemma 2.2, one can choose nonnegative
functions 
$\theta^-\in\cic(H^-)$, $\theta^0\in\cic(H^0)$ 
with 
$$
\int_{H^-}  \theta^-\,
 =
\int_{H^0} \theta^0\, 
= 1\tag 2.5 
$$ such that
$$
\text{supp}(\theta^-)\cdot\text{supp}
(\theta^0)\subset B_{\tilde H}({r})\,,\tag 2.6 
$$
and at the same time
$$
\|\tilde \theta\|_{\ell}\ll
{r}^{(2\ell+N/2)}\,,\tag 2.7
$$
where $\tilde \theta\in\cic(\tilde H)$ is defined by $$\tilde
\theta(h^-h^0)\df \theta^-(h^-)\theta^0(h^0)\Delta(h^0)^{-1}\,.
$$  Also define $\varphi\in \cic(X)$ by 
$\varphi(h^-h^0hz) = \tilde\theta(h^-h^0)f(h)$; the definition makes sense
because of (2.6) and assumptions (i), (ii) of the theorem. 
Then $I_{f,\psi}(g_t,z)$ can be reasonably well approximated by $\langle
g_t\varphi,\psi\rangle = \langle \varphi,g_{-t}\psi\rangle$: 
$$
\aligned
&\left|I_{f,\psi}(g_t,z) - \langle \varphi,g_{-t}\psi\rangle\right|\\
&\un{(2.3)}=\left|\int_H f(h)\psi(g_thx)\, d\nu (h) - \int_G \tilde
\theta(h^-h^0)
f(h) \psi(g_th^-h^0hx)\,d\mu(h^-h^0h)\right|\\
&\un{(2.5)}=\left|\int_G \theta^-(h^-) \theta^0(h^0)
f(h) \Big(\psi(g_thx) -
\psi\big(\Phi_t(h^-h^0)g_thx\big)\Big)
\Delta(h^0)^{-1}\,d\mu(h^-h^0h)\right|\\  &\un{(2.2),\,(2.6) }\le
\sup_{g\in B_{\tilde H}(r),\,y\in X}|\psi(gy) -
\psi(y)|\int_G \left|\theta^-(h^-) \theta^0(h^0)
f(h)
\Delta(h^0)^{-1}\right|\,d\mu(h^-h^0h) \\  &\un{(2.3)}\le  
\|\psi\|_{\Lip}\cdot r \cdot\int_H|f|\,
\,. \endaligned
$$
On the other hand, in view of
Lemma 2.2 and $\pi_z$ being a local isometry,
$$
\|\varphi\|_{\ell} = \|\tilde \theta\cdot f\|_{\ell}
\ll\|\tilde
\theta\|_{\ell}\|f\|_{\ell} 
\un{(2.7)}\ll {r}^{-(2\ell+N/2)} \|f\|_{\ell}\,, 
$$ 
hence  by (2.1) 
$$
\left|\langle g_t\varphi,\psi\rangle 
\right| 
\ll {r}^{(2\ell+N/2)} \|f\|_{\ell}\|\psi\|_{\ell}e^{-\gamma t}\,, 
$$
finishing the proof. \qed\enddemo

\example{Remark 2.4} In order to derive Theorem 1.1 from Theorem 2.3 
it suffices to choose $r = e^{-\beta t}$ for some suitable $\beta$.
The same  trick will help us in the proof of Theorem 1.3.
Note that  $t$ needs to be taken large enough so that condition (ii) of
Theorem 2.3 is satisfied for all $z\in L$. The latter is possible
because, in view of the compactness of $L$ and discreteness of
$\Gamma$ in
$G$, the value
$$r(L)\df\inf_{z\in L}\sup\{r > 0\mid \pi_z: G\to X\text{ is injective on
}B(r)\}$$
is positive; we will call it the {\sl injectivity
radius of\/} $L$.
\endexample

\example{Remark 2.5} It is worthwhile to point out that 
 $H$ being the \ehs\ relative to $g_t$, $t > 0$, 
was crucially important for the proof. When $g_t$
is replaced with $g_\vt$ where $\vt$ is an arbitrary element of $\fa^+$, 
one can still talk about 
$\Phi_\vt$, the inner automorphism of $H$
given by  $$\Phi_\vt(h)\df g_\vt h
(g_\vt)^{-1}\,.\tag 2.8
$$ It is expanding on $H$, since the 
latter is contained in the \ehs\ relative to $g_\vt$;
however it is not  non-expanding on $\tilde H$ in the sense of (2.2), 
thus there is no guarantee that  $I_{f,\psi}(g_t,z)$ is close to
$ \langle \varphi,g_{-t}\psi\rangle$ for $\varphi$ constructed as in the 
above proof. We bypass this difficulty by means of 
an additional step, based on the nondivergence phenomenon, to be described
in the next section.
\endexample


\comment

\subheading{2.5. Localizing by smooth test functions} Finally, since we
need to apply Theorem 2.1 to functions with small support, the following
elementary lemma will be useful:

 \proclaim{Lemma 2.7 \cite{KM1, Lemma 2.4.7(b)}} 
For any
 $0<r<1$, 
there exists a nonnegative function $\theta\in
\cic(H)$ such that 
{\rm supp}$(\theta)\subset B(e,r)$, 
$\int_H
\theta\,d\nu = 1$, and
$$\|\theta\|_{\ell}\le \tilde c  r^{-(\ell+mn/2)}\,,$$ where $\tilde c$
is a constant independent on $r$.
\endproclaim

 We will fix a
right-invariant metric `dist' on $G$, giving rise to the corresponding
metric on $X$.
$B(x,r)$ will stand for an open ball of radius $r$ centered at $x$.
For any
$z\in X$ we define the {\sl injectivity radius\/} $r(z)$ at $z$ to be
the supremum of
$r' > 0$ for which the map $\pi_z: G\to X$ is injective on $B(e,r')$. It
follows from
the discreteness of $\Gamma$ in $G$ that $r$ is a proper function,
that is, for any compact $L\subset X$ the value $r(L)
\df\inf_{z\in L}r(z)$ is positive; we will call it the {\sl injectivity
radius of\/} $L$.

\comdima{It seems to me that this is what can be proved using 
our argument, do we  need to write it down? Perhaps it will be better if
we do..}  Here $\|\cdot\|_{\ell}$ stands for the norm in the Sobolev space
$W_{\ell}^2(H,\nu)$. It is possible to explicitly estimate $E$ as
well (through $\|\psi\|_{\ell}$ as well as other parameters)
but it is not needed for the main result of the paper.

Consequently, a uniform  estimate with $z$ in some compact set 
will involve the injectivity radius of a slightly bigger compact set in
some negative power.
\endcomment

\heading{3.  Quantitative nondivergence}\endheading

For any $\vre > 0$ consider
$$
K_{{\vre}} \df \pi\big(\big\{ g \in G \bigm| \| g
\vv \| \ge  {\varepsilon} \quad \forall\, \vv \in
\bz^{k}\nz\big\}\big)\,.
$$
In other words, $
K_{{\vre}}$ consists of lattices in $\br^k$ with no nonzero vector
of length less than $\vre$. These sets are compact by virtue of Mahler's
Compactness Criterion (see \cite{R,~Corollary 10.9} or \cite{BM}). 
Here $\|\cdot\|$ can be any norm on $\br^k$ which we will
from now on take to be the standard Euclidean norm.

It was proved in \cite{KM2}, refining previous work on non-divergence
of unipotent flows
\cite{M, D}, that certain polynomial maps from 
balls in Euclidean 
spaces to $X$ cannot take values outside of $K_\vre$ on a set of big
measure. Namely, the following is a 
 special case of 
\cite{BKM, Theorem 6.2} (see also \cite{KLW, KT, K} for further
generalizations):

 \proclaim{Theorem 3.1} For $d\in\bn$,  let $\varphi$ be a map 
 $\br^d \to \GL_k(\br)$ 
such that 
\roster
\item"(i)" all coordinates (matrix elements) of $\varphi(\cdot)$
 are affine (degree
$1$ polynomials), 
\endroster
and let  a ball $B\subset\br^d$
and  
$0 < \rho  \le 1$ be such that 
\roster
\item"(ii)" for any $j = 1,\dots,k-1$ and any 
$\vw \in
\bigwedge^j(\bz^{k})\nz$ one has
$$
\| \varphi(x)
\vw \| \ge  \rho\quad\text{for some }x\in B\,.
$$
\endroster
Then 
 for any  positive $ \vre \le \rho$ one has
$$
\lambda\left(\big\{x\in B\mid \pi\big(\varphi(x)\big)\notin K_\vre
\big\}\right)
\ll
 \left(\frac\vre \rho \right)^{1/d(k-1)}  \lambda(B)\,.
\tag 3.1
$$
\endproclaim
 
Here $\lambda$ is Lebesgue measure on $\br^d$, and the
Euclidean\footnote{In \cite{KM2} the statement of Theorem 5.2 involved 
the sup norm instead of the Euclidean one, which resulted in a
restriction
for $\rho$ to be not greater than $1/k$; thus we chose to refer to 
\cite{BKM} for the Euclidean norm version.} norm
$\|\cdot\|$ is naturally extended from
$\br^k$ to its exterior powers. We remark that  the way 
assumption (i)
is used in the proof is by verifying  that 
 all the functions
$x\mapsto
\|
\varphi(x)
\vw \|$, where $\vw \in
\bigwedge^j(\bz^{k})$, are \cag\ on $\br^d$,
 with  some fixed $C = C(d,k) > 0$ and $\alpha =
1/d(k-1)
$, the exponent appearing in (3.1). 
See
\cite{KM2} for more detail.


\medskip

Our plan is  to apply Theorem 3.1 with $\varphi:\mr\to G$ given by
$$
\varphi(Y) = g_\vt u_{ \sssize Y}g\tag 3.2
$$ for some  $g\in G$ and $\vt\in\fa^+$.
It is clear that assumption (i) holds. As for 
(ii),
we will need to have uniformity in $\vt\in\fa^+$ and in $g$ such that
$\pi(g)$ belongs to a compact subset of $X$. This can be 
extracted from the next lemma, which is
immediate from \cite{KW, 
Proposition 2.4} applied
to  the representations of $G$ on $\bigwedge^j(\br^{k})$, $\,j =
1,\dots,k-1$:

 \proclaim{Lemma 3.2} There exists $\alpha > 0$ with the following
property. Let $B$ be a ball centered at 
$0$ in $\mr$. Then one can find $b 
> 0$ such that 
 for any $\,j = 1,\dots,k-1$,  any 
$\vw \in
\bigwedge^j(\br^{k})$ 
 and any $\,\vt\in\fa^+$ 
one has 
$$
\sup_{Y \in B} \big\| g_{\vt}u_{ \sssize Y}\vw \big\| \ge 
be^{\alpha\lfloor \vt\rfloor}\|\vw\|\,.
$$
\endproclaim


 \proclaim{Corollary 3.3} 
Let $B$ be a neighborhood of 
$0$ in $\mr$ and let $L\subset X$ be compact. Then there exists
$b > 0$ such that 
 for any $j = 1,\dots,k-1$,  any 
$\vw \in
\bigwedge^j(\bz^{k})\nz$, any $g\in\pi^{-1}(L)$ 
 and any $\,\vt\in\fa^+$ 
one has 
$$
\sup_{Y \in B} \big\| g_{\vt}u_{ \sssize Y}g\vw \big\| \ge 
b e^{\alpha\lfloor \vt\rfloor}\,.
$$
\endproclaim

\demo{Proof} Apply the above lemma with $\vw$ replaced by $g\vw$; 
it follows from the compactness of $L$ and discreteness of
$\bigwedge^j(\bz^{k})$
in $\bigwedge^j(\br^{k})$ that
$$\inf\big\{\|g\vw\|
\bigm|
\pi(g)\in L,\ \vw \in
\tsize\bigwedge^j(\bz^{k})\nz, j = 1,\dots,k-1\big\} $$ is positive.
\qed\enddemo

 \proclaim{Corollary 3.4} 
Let $L\subset X$ be compact and let $B\subset H$ be a ball centered at 
$e\subset H$.  Then there exists 
$T = T(B,L)$ such that for every
$0 < \vre < 1$, any $z\in L$ and any
$\,\vt\in\fa^+$ with $\lfloor \vt\rfloor \ge T$ one has
$$
\nu\left(\big\{h\in B\mid g_{\vt}h z\notin K_{
\vre}
\big\}\right)
\ll 
\vre^{\frac{1}{mn(k-1)}} \nu(B)\,.
$$
\endproclaim


\demo{Proof} 
Define $T$  by
$b e^{\alpha T} = 1$, where $\alpha$ is given by Lemma 3.2 
and $b$ by Corollary
3.3 applied to $\log(B)\subset \mr$ and $L$. (Note that the exponential
map from $\mr$ to $H$ is an isometry.)   Take
$\varphi$  as in (3.2) with
$g\in\pi^{-1}(L)$. Then, in view of Corollary 3.3, assumption (ii) of
Theorem 3.1, with $d = mn$, will be satisfied with
$\rho = 1$ as long as $\lfloor \vt\rfloor \ge T$. 
\qed\enddemo

We conclude this section by an estimate of 
the injectivity radius of $K_\vre$, to make it possible to combine the
above  corollary with Theorem 2.3.
Observe that any lattice $\Lambda\in K_\vre$ can be generated by
vectors of norm
$\ll 1/\vre^{k-1}$;  if $g\Lambda = \Lambda$ and $g\ne e$, then for one of
those vectors
$\vv$ one has $\|g\vv - \vv\| \ge \vre$. This implies that dist$(e,g) \gg
\vre^k$. 
We
arrive at

 \proclaim{Proposition 3.5} There exists positive $c = c(k)$ such that
$r(K_\vre)
\ge c
\cdot
\vre^k$ 
$\forall\,\vre>0$.
\endproclaim

\heading{4.  Proof of Theorem 1.3}\endheading

Our goal in this section will be to find $\tilde \gamma > 0$ such that
for any $f\in \cic(H)$, $\psi\in
W^{2,\infty}(X)\cap \Lip(X)$ with $\int_X\psi = 0$ and   compact  $L\subset X$
there exists $\tilde C> 0$ such that for all $z\in L$ and
all $\vt\in\fa^+$ one has
$$
\left|I_{f,\psi}(g_\vt,z) \right|\le \tilde C e^{-\tilde\gamma
\lfloor
 \vt\rfloor}\,.\tag 4.1  
$$
Then Theorem 1.3 will follow by applying (4.1) with $\psi$ replaced by
$\psi - \int_X\psi$. 
 \comment
Since both 
 supp$(f)$ and $L$ are
compact, 
$f$ can be
written as a sum of 
functions $f_j$, $1\le j\le N$, with $\pi_z$ injective on supp$(f_j)$ for
all $z\in L$ and for each $j$. Hence one can without loss of
generality (by increasing $\tilde C$ if needed) assume that  supp$(f)$ 
is contained in a ball $B = B(e,\rho)$  
and maps $\pi_z$ are injective on $B$ for all
$z\in L$. 
\endcomment
Note also that,  by increasing $\tilde
C$ if needed, it is enough to prove (4.1) for $\vt$ with large
enough  
$\lfloor\vt\rfloor$.

Given $\vt\in\fa^+$, define $t \df \lfloor\vt\rfloor/2$, 
and let $$\vu =  \vu(\vt) \df \vt -
\big(\tfrac{t}{m},\dots,\tfrac{t}{m},  
\tfrac{t}{n},\dots,\tfrac{t}
{n}\big)\,.\tag 4.2
$$ Note that
$\vu\in\fa^+$,
$\lfloor\vu\rfloor \ge \lfloor\vt\rfloor/2 = t$, and $g_\vt = g_t  g_\vu$
(here $g_\vt$ and $  g_\vu$ are defined via (1.4), and $g_t$ is as in
(1.2)).

\comment Using the fact that  $g_\vt = g_t  g_\vu$, we organize the proof
in two steps.
\comdima{what follows is a little user's guide to the proof, feel free to
edit or delete it}
First we apply $g_\vu$ and use 
 Corollary 3.4 to conclude that most of the translate
$g_\vu Bz$ is contained in some fixed compact set. Then we divide 
$g_\vu Bz$ into little pieces  and use 
the equidistribution 
of the $g_t$-translates of most of the pieces, that is, Theorem 2.1.
\endcomment

Take a function $\theta$ supported on $B_H(r)$ as in Lemma 2.2(a), with
$\,r = e^{-\beta t}$ where
$\beta$ is to be specified later; since $\int_H
\theta
= 1$ and $\nu$ is translation-invariant, one can 
 write
$$
\aligned
I_{f,\psi}(g_\vt,z) &= 
\int_H
f(h)\psi(g_\vt hz)\,d\nu(h)
\int_H
\theta(y)\,d\nu(y)\\&= 
\int_H
\int_H
f\big(\Phi_\vu^{-1}(y)h\big)\theta(y)\psi\big(g_t  g_\vu
\Phi_\vu^{-1}(y)hz\big)\,d\nu(y)\,d\nu(h)\\&=\int_H
\int_H
f\big(\Phi_\vu^{-1}(y)h\big)\theta(y)\psi\big(g_t y 
g_\vu hz\big)\,d\nu(y)\,d\nu(h)\,.
\endaligned
 $$
Note that $\Phi_\vu^{-1}$ is a
contracting
automorphism of $H$, in fact, one has $$\dist\big(e,\Phi_\vu^{-1}(h)\big)
\le e^{-2\lfloor\vu\rfloor}\dist(e,h) \le e^{-2t}\dist(e,h)
$$
for any $h\in H$.
Choose $B = B(r)$ containing $\supp f$. Then
  the supports of all  functions of the form $h\mapsto
f\big(\Phi_\vu^{-1}(y)h\big)$ are contained in  $$\tilde B \df B(r +
e^{-(2+\beta)t})\,.$$ 
By taking $t$ large enough it is safe to assume that 
$$
e^{-\beta t}  < r_0/2\,,\tag 4.3$$
$\nu(\tilde B)\le 2\nu(B)$, and also that 
$t > T \df T(\tilde B,L)$ as in  Corollary 3.4. Now define 
$\vre$ by
$$
\vre = \left(\tfrac {2}{c}e^{-\beta t} \right)^{1/k}\,,\tag 4.4
$$
where $c$ is from Proposition 3.5, and denote
$$A \df \big\{h\in \tilde B\mid g_{\vu}h z\notin K_{
\vre}
\big\}\,.$$ Then for any $\vu\in\fa^+$ with $\lfloor\vu\rfloor \ge T$
 and any $z\in L$ one knows, in view of  Corollary 3.4, that 
$$\nu(A) \ll 
\vre^{\frac{1}{mn(k-1)}} \nu(\tilde B)\,.
$$
Hence the absolute value of 
$$
\int_{A}\int_H
f\big(\Phi_\vu^{-1}(y)h\big)\theta(y)\psi\big(g_t y 
g_\vu hz\big)\,d\nu(y)\,d\nu(h)
$$
is 
$$ 
\ll 
\vre^{\frac{1}{mn(k-1)}} \nu(\tilde B)\sup|f|\sup|\psi|\int_H
\theta\,
\ll 
\sup|f|\sup|\psi| \nu(B)\cdot e^{-\frac{{\beta }}{mnk(k-1)}t}\,.
 $$

Now let us assume that $h\in \tilde B\ssm A$, and apply Theorem 2.3
with $r = e^{-\beta t}$, $
g_\vu hz $ 
in place of $z$ and $$f_h(y)
\df  f\big(\Phi_\vu^{-1}(y)h\big)\theta(y)$$ in place of $f$.
Clearly condition (i) 
follows from (4.3), and, since  $g_\vu hz  \in K_\vre$ whenever 
$h\notin A$, condition (ii) is satisfied in view of 
Proposition 3.5 and (4.4). 
Also, because
$\Phi_\vu^{-1}|_H$ is contracting,   partial derivatives of $f$
will not increase after precomposition with $\Phi_\vu^{-1}$, and thus
$$\|f_h\|_\ell \un{Lemma 2.2(b)}\ll \| f\|_{C^\ell} \|\theta\|_\ell 
\un{Lemma 2.2(a)}\ll   
r^{-(\ell+mn/2)}\| f\|_{C^\ell}\,.\tag 4.5$$
This way one gets:
$$
\aligned
\left|\int_{\tilde B \ssm A}\int_H
f\big(\Phi_\vu^{-1}(y)h\big)\theta(y)\psi\big(g_t y 
g_\vu hz\big)\,d\nu(y)\,d\nu(h)\right| \le 
\int_{\tilde B \ssm A}\left|I_{f_h,\psi}(g_t,g_\vu h z)
\right|\,d\nu(h)\\ \un{(2.4)}\le  
 E(\psi) \left(r 
\int_H |f_h|\,
 +
r^{-(2\ell + N/2)} \|f_h\|_{\ell}
e^{-\gamma
t}\right) \nu(\tilde B)
\\ \un{(4.5)}\ll   E(\psi)
\left(\sup|f|\cdot e^{-\beta t} 
 +
\| f\|_{C^\ell}  \cdot 
e^{-(\gamma - (2\ell + N/2)\beta)
t}\right) \nu( B)\,.\endaligned$$ 

Combaning the two estimates above, one can conclude that
$$
\split
\left|
I_{f,\psi}(g_\vt,z)\right|& \ll C_1 e^{-\frac{{\beta
}}{mnk(k-1)}t} +  C_2 e^{-\beta t} 
 +
C_3 
e^{-(\gamma - (2\ell + N/2)\beta)
t}\\
& \le \max(C_1,C_2) e^{-\frac{{\beta
}}{mnk(k-1)}t} 
 +
C_3 
e^{-(\gamma - (2\ell + N/2)\beta)
t}
\,,
\endsplit$$
where $C_i$, $i = 1,2,3$, depend on $f$, $\psi$ and $L$. 
An elementary computation shows that  choosing $\beta$  equalizing the
two exponents
above will produce $$
\tilde \gamma = \frac \gamma {1 + mnk(k-1)(2\ell+N/2)}
 $$ such that  (4.1) will hold with  $\tilde C \ll \max(C_1,C_2,C_3)$.
\qed
\comment
$$
\big(\const_1 r^{\frac{1}{mnk(k-1)}} + \const_2 
 r^{-(s+\ell+mn/2)}e^{-\gamma t}\big)\nu(B)\,,
$$
where the constants above depend on $f$ and $\psi$. 
An elementary computation shows that if one chooses $r = e^{-\beta t}$
where
$\beta$ is such that the exponents of the two summands above are the
same,  one would get (1.5) with some $\tilde C$ and 
$$
\tilde \gamma = \frac \gamma {1 + mnk(k-1)(s+\ell+mn/2)}\,. 
 $$
\endcomment

\enddocument